\documentclass[10pt,twoside,reqno]{amsart}
\usepackage{amssymb}
\usepackage{multirow}
\calclayout

\numberwithin{equation}{section}
\makeatletter

\renewcommand{\@secnumfont}{\bfseries}

\renewcommand{\section}{\@startsection{section}{1}%
  {0mm}{.7\linespacing\@plus\linespacing}{.5\linespacing}
  {\normalfont\bfseries\centering}}

\newcommand{\bibsection}{\@startsection{section}{1}%
  {0mm}{.7\linespacing\@plus\linespacing}{.5\linespacing}
  {\normalfont\scshape\centering}}

\renewcommand{\@biblabel}[1]{#1.}

\newtheorem{thm}{\bf Theorem}[section]

\begin{document}

\vspace{1.3cm}

\title[Extended Stirling and extended Bell polynomials]{Extended Stirling polynomials of the second kind and extended Bell polynomials}

\author{Taekyun Kim}
\address{Department of Mathematics, Kwangwoon University, Seoul 139-701, Republic
	of Korea}
\email{tkkim@kw.ac.kr}

\author{Dae San Kim}
\address{Department of Mathematics, Sogang University, Seoul 121-742, Republic of Korea}
\email{dskim@sogang.ac.kr}


\begin{abstract} 
Recently, several authors have studied the Stirling numbers of the second kind and Bell polynomials. In this paper, we study the extended Stirling polynomials of the second kind and the extended Bell polynomials associated with the Stirling numbers of the second kind. In addition, we note that the extended Bell polynomials can be expressed in terms of the moments of the Poisson random variable with parameter $\lambda >0$.
\end{abstract}
\maketitle

\section{Introduction}
As is well known, the Stirling numbers of the second kind are defined as
\begin{equation}\begin{split}\label{01}
x^n = \sum_{l=0}^n S_2(n,l) (x)_l\,\, (n \geq 0), \quad (\textnormal{see} \,\, [1-16]).
\end{split}\end{equation}
The generating function of $S_2(n,l)$ is given by
\begin{equation}\begin{split}\label{02}
\frac{1}{m!} (e^t-1)^m = \sum_{n=m}^\infty S_2(n,m) \frac{t^n}{n!},
\end{split}\end{equation}
where $m \in \mathbb{N} \cup \{0\}$, (see [2,7,8]).

The Stirling polynomials of the second kind are defined by the generating function 
\begin{equation}\begin{split}\label{03}
\frac{1}{k!} e^{xt} (e^t-1)^k = \sum_{n=k}^\infty S_2(n,k|x) \frac{t^n}{n!},
\end{split}\end{equation}
where $k \geq 0$, (see [3,5,14]).

From \eqref{02} and \eqref{03}, we note that
\begin{equation}\begin{split}\label{04}
S_2(n,k|x) &= \sum_{l=k}^n {n \choose l} S_2(l,k) x^{n-l}\\
&= \sum_{l=0}^{n-k} {n \choose l} S_2(n-l,k) x^l,
\end{split}\end{equation}
where $n, k \geq 0$, (see [3,4,5,14]).

The Bell polynomials are defined by the generating function 
\begin{equation}\begin{split}\label{05}
e^{x(e^t-1)} = \sum_{n=0}^\infty Bel_n(x) \frac{t^n}{n!},\quad (\textnormal{see} \,\, [7,8,9]).
\end{split}\end{equation}
When $x=1$, $Bel_n(1) = Bel_n$, $(n \geq 0)$, are called the Bell numbers.

From \eqref{02} and \eqref{05}, we note that

\begin{equation}\begin{split}\label{06}
e^{x(e^t-1)} &= \sum_{m=0}^\infty x^m \frac{1}{m!} (e^t-1)^m\\
&= \sum_{n=0}^\infty \left( \sum_{m=0}^n S_2(n,m) x^m \right) \frac{t^n}{n!}.
\end{split}\end{equation}
Thus, by \eqref{06}, we get
\begin{equation*}\begin{split}
Bel_n(x) =  \sum_{m=0}^n S_2(n,m) x^m ,\,\,(n \geq 0).
\end{split}\end{equation*}
A random variable $X$, taking on one of the values 0,1,2,$\cdots$, is said to be a Poisson random variable with parameter $\lambda >0$ if $P(i) = P(X=i) = e^{-\lambda} \frac{\lambda^i}{i!}$, $i=0,1,2,\cdots$. Note that $\sum_{i=0}^\infty P(i) = e^{-\lambda} \sum_{i=0}^\infty \frac{\lambda^i}{i!} = e^{-\lambda}\cdot e^{\lambda} =1$.

The expectation of a Poisson random variable with parameter $\lambda$ is given by
\begin{equation}\begin{split}\label{07}
E[X] = \sum_{i=0}^\infty iP(i) = \sum_{i=0}^\infty i e^{-\lambda} \frac{\lambda^i}{i!} = \lambda.
\end{split}\end{equation}
The moments of Poisson random variable $X$ with parameter $\lambda >0$ is defined by 
\begin{equation}\begin{split}\label{08}
E[X^n]= \sum_{x=0}^\infty x^n P(x) = e^{-\lambda } \sum_{x=0}^\infty x^n \frac{\lambda ^x}{x!},
\end{split}\end{equation}
where $n \in \mathbb{N}$ (see [15]).

When $n=1$, the first moment $E[X]$ is the mean (or expection) of $X$ with parameter $\lambda >0$. Recently, several authors have studied the Stirling numbers of the second kind and Bell polynomials (see [5-16]).

In this paper, we consider the extended Stirling polynomials of the second kind and the extended Bell polynomials associated with the Stirling numbers of the second kind. Then we give some identities between the extended Stirling numbers of the second kind and the extended Bell polynomials. From our new identities and properties of those numbers and polynomials, we note that the extended Bell polynomials can be expressed in terms of the moments of the Poisson random variable with parameter $\lambda >0$.

\section{Extended Stirling polynomials of the second kind and extended Bell polynomials}

For $k >0$, we define the extended Stirling polynomials of the second kind given by the generating function 
\begin{equation}\begin{split}\label{09}
\frac{1}{k!} e^{xt} (e^t-1+rt)^k = \sum_{n=k}^\infty S_{2,r}(n,k|x) \frac{t^n}{n!},
\end{split}\end{equation}
where $x,r \in \mathbb{R}$.

When $x=0$, $S_{2,r}(n,k|0) = S_{2,r}(n,k)$, $(n,k \geq 0)$, are called the extended Stirling numbers of the second kind. Note that $S_{2,0}(n,k)=S_2(n,k)$ are the Stirling numbers of the second kind. 

It is easy to show that
\begin{equation}\begin{split}\label{10}
\frac{1}{k!} e^{xt} (e^t-1+rt)^k = \sum_{n=k}^\infty \left( \sum_{m=k}^n {n \choose m} S_{2,r}(m,k) x^{n-m} \right) \frac{t^n}{n!}.
\end{split}\end{equation}
By \eqref{09} and \eqref{10}, we get
\begin{equation}\begin{split}\label{11}
S_{2,r}(n,k|x) = \sum_{m=k}^n {n \choose m} S_{2,r}(m,k) x^{n-m},
\end{split}\end{equation}
where $n,k \geq 0$ and $r,x \in \mathbb{R}$.

We observe that
\begin{equation}\begin{split}\label{12}
\sum_{k=0}^\infty \frac{1}{k!} (e^t-1+rt)^k = \sum_{n=0}^\infty \left( \sum_{k=0}^n S_{2,r}(n,k) \right) \frac{t^n}{n!}, 
\end{split}\end{equation}
and
\begin{equation}\begin{split}\label{13}
\sum_{k=0}^\infty \frac{1}{k!} (e^t-1+rt)^k = e^{e^t-1+rt}.
\end{split}\end{equation}
In view of \eqref{05}, we can define the extended Bell numbers which are given by the generating function
\begin{equation}\begin{split}\label{14}
e^{e^t-1+rt} = \sum_{n=0}^\infty Bel_{n,r} \frac{t^n}{n!},
\end{split}\end{equation}
From \eqref{12} and \eqref{14}, we have
\begin{equation}\begin{split}\label{15}
Bel_{n,r} = \sum_{k=0}^n S_{2,r} (n,k), \,\,(n \geq 0).
\end{split}\end{equation}
Note that $Bel_{n,0}=\sum_{k=0}^n S_{2,0}(n,k) = \sum_{k=0}^n S_2(n,k)$. Now, we define the extended Bell polynomials given by the generating function as follows:

\begin{equation}\begin{split}\label{16}
e^{\lambda (e^t-1+rt)} = \sum_{n=0}^\infty Bel_{n,r}(\lambda ) \frac{t^n}{n!},
\end{split}\end{equation}
where $\lambda ,r\in \mathbb{R}$.

From \eqref{16}, we note that
\begin{equation}\begin{split}\label{17}
e^{\lambda (e^t-1+rt)} &= \sum_{m=0}^\infty \lambda ^m \frac{1}{m!} (e^t-1+rt)^m\\
&= \sum_{m=0}^\infty \lambda ^m \sum_{n=m}^\infty S_{2,r}(n,m) \frac{t^n}{n!}\\
&= \sum_{n=0}^\infty \left( \sum_{m=0}^n \lambda ^m S_{2,r}(n,m) \right) \frac{t^n}{n!}.
\end{split}\end{equation}
Therefore, we obtain the following theorem.

\begin{thm}
For $n \geq 0$, we have
\begin{equation*}\begin{split}
Bel_{n,r}(\lambda ) = \sum_{m=0}^n \lambda ^m S_{2,r}(n,m),
\end{split}\end{equation*}
and
\begin{equation*}\begin{split}
S_{2,r}(n,m|x) = \sum_{k=m}^n  {n \choose k} S_{2,r}(k,m) x^{n-k},
\end{split}\end{equation*}
where $n, m \geq 0$ and $r \in \mathbb{R}$.
\end{thm}

From \eqref{01}, we note that

\begin{equation}\begin{split}\label{18}
\sum_{n=k}^\infty S_{2,r}(n,k) \frac{t^n}{n!}&=  \frac{1}{k!} ( e^t-1+rt)^k = \frac{1}{k!}\sum_{l=0}^k {k \choose l} r^l t^l (e^t-1)^{k-l}\\
&=\frac{1}{k!} \sum_{l=0}^k \frac{k!}{l!(k-l)!} r^l t^l (e^t-1)^{k-l}\\
&= \sum_{l=0}^k \frac{r^l}{l!}t^l \sum_{n=k}^\infty S_2(n-l,k-l) \frac{t^{n-l}}{(n-l)!}\\
&= \sum_{n=k}^\infty \left( \sum_{l=0}^k {n \choose l} r^l S_2(n-l,k-l) \right) \frac{t^n}{n!}.
\end{split}\end{equation}

Thus, by comparing the coefficients on both sides of \eqref{18}, we obtain the following theorem.

\begin{thm}
For $n,k \geq 0$, we have
\begin{equation*}\begin{split}
S_{2,r}(n,k) = \sum_{l=0}^k {n \choose l} r^l S_2(n-l,k-l).
\end{split}\end{equation*}
\end{thm}

It is not difficult to show that

\begin{equation}\begin{split}\label{19}
e^{e^t-1+rt} = \sum_{n=0}^
\infty \left( \sum_{l=0}^n {n \choose l} Bel_l r^{n-l} \right) \frac{t^n}{n!}.
\end{split}\end{equation}

Thus, by \eqref{14} and \eqref{19}, we easily get
\begin{equation}\begin{split}\label{20}
Bel_{n,r} = \sum_{k=0}^n S_{2,r} (n,k) = \sum_{k=0}^n {n \choose k} Bel_k r^{n-k},
\end{split}\end{equation}
where $r \in \mathbb{R}$ and $n \in \mathbb{N} \cup \{0\}$.

From \eqref{20}, we have
\begin{equation*}\begin{split}
Bel_{n,r} = \sum_{l=0}^n {n \choose l} Bel_{n-l} r^l = \sum_{k=0}^n \left( \sum_{l=0}^k {n \choose l} r^l S_2(n-l,k-l) \right).
\end{split}\end{equation*}
By \eqref{02}, we get
\begin{equation}\begin{split}\label{22}
&\sum_{n=m}^\infty S_2(n,m) \frac{t^n}{n!} = \frac{1}{m!} (e^t-1)^m = \frac{1}{m!} (e^t -1+rt-rt)^m\\
&= \frac{1}{m!} \sum_{l=0}^m {m \choose l} t^l (e^t-1+rt)^{m-l} (-1)^l r^l\\
&= \frac{1}{m!} \sum_{l=0}^m \frac{m! t^l (-1)^l r^l}{l! (m-l)!} (e^t-1+rt)^{m-l}\\
&= \sum_{l=0}^m  \frac{r^l (-1)^l}{l!} t^l \sum_{n=m-l}^\infty S_{2,r}(n,m-l) \frac{t^n}{n!}\\
&= \sum_{n=m}^\infty \left( \sum_{l=0}^m S_{2,r} (n-l,m-l) {n \choose l} (-1)^l r^l \right) \frac{t^n}{n!},
\end{split}\end{equation}
where $m \in \mathbb{N}\cup \{0\}$.
By comparing the coefficients on both sides of \eqref{22}, we obtain the following theorem.

\begin{thm}
For $n \geq m \geq 0$ and $r \in \mathbb{R}$, we have
\begin{equation*}\begin{split}
S_2(n,m)=\sum_{l=0}^m {n \choose l} (-1)^l r^l S_{2,r} (n-l,m-l) 
\end{split}\end{equation*}
\end{thm}

Now, we consider the inversion formula of \eqref{22}. From \eqref{09}, we note that
\begin{equation}\begin{split}\label{23}
	&\frac{1}{k!} (e^t -1+rt)^k = \frac{1}{k!} \sum_{l=0}^k {k \choose l} r^l t^l (e^t-1)^{k-l}\\
	&= \frac{1}{k!} \sum_{l=0}^k r^l t^l \frac{k!}{l!(k-l)!} (e^t-1)^{k-l}\\
	&= \sum_{l=0}^k  \frac{r^l}{l!} t^l \sum_{n=k-l}^\infty S_2(n,k-l) \frac{t^n}{n!}\\
	&= \sum_{n=k}^\infty \left( \sum_{l=0}^k {n \choose l} r^l S_2(n-l,k-l) \right) \frac{t^n}{n!}.
	\end{split}\end{equation}	
Therefore, by \eqref{09} and \eqref{23}, we obtain the following theorem.

\begin{thm}
For $n\geq k \geq 0$, we have
\begin{equation*}\begin{split}
S_{2,r}(n,k) =  \sum_{l=0}^k {n \choose l} r^l S_2(n-l,k-l) .
\end{split}\end{equation*}
\end{thm}
For $m,n,k \in \mathbb{N}$, we have
\begin{equation}\begin{split}\label{24}
&\frac{1}{m!} (e^t-1+rt)^m \frac{1}{k!} (e^t-1+rt)^k = \frac{1}{m!k!} (e^t-1+rt)^{m+k}\\
&= \frac{(m+k)!}{m!k!} \frac{1}{(m+k)!} (e^t-1+rt)^{m+k} = {m+k \choose m} \sum_{n=m+k}^\infty S_{2,r} (n,m+k) \frac{t^n}{n!}.
\end{split}\end{equation}
On the other hand,
\begin{equation}\begin{split}\label{25}
&\frac{1}{m!} (e^t-1+rt)^m \frac{1}{k!} (e^t-1+rt)^k\\
&= \left( \sum_{l=m}^\infty S_{2,r}(l,m) \frac{t^l}{l!} \right) \left( \sum_{j=k}^\infty S_{2,r}(j,k) \frac{t^j}{j!} \right)\\
&= \sum_{n=k+m}^\infty \left( \sum_{l=m}^n \frac{n!}{l! (n-l)!}	S_{2,r}(l,m) S_{2,r}(n-l,k) \right) \frac{t^n}{n!}\\
&= \sum_{n=k+m}^\infty \left( \sum_{l=m}^n {n \choose l} S_{2,r}(l,m) S_{2,r}(n-l,k) \right) \frac{t^n}{n!}.
\end{split}\end{equation}
Therefore, by \eqref{24} and \eqref{25}, we obtain the following theorem.

\begin{thm}
For $n,m,k \geq 0$ with $n \geq m+k$, we have
\begin{equation*}\begin{split}
{m+k \choose m} S_{2,r}(n,m+k)= \sum_{l=m}^n {n \choose l} S_{2,r}(l,m) S_{2,r}(n-l,k).
\end{split}\end{equation*}
\end{thm}

Now, we observe that
\begin{equation*}\begin{split}
&\frac{1}{m!} (e^t-1+rt)^m \frac{1}{k!} (e^t-1+rt)^k\\
&=\left( \frac{1}{m!} \sum_{l=0}^m {m \choose l} (e^t-1)^{m-l} r^l t^l \right) 
\left( \frac{1}{k!} \sum_{j=0}^k {k \choose j} (e^t-1)^{k-j} r^j t^j \right) \\
&= \left( \sum_{l=0}^m \frac{r^l t^l}{l!} \sum_{n_1=m}^\infty S_2(n_1-l,m-l) \frac{t^{n_1-l}}{(n_1-l)!} \right)\\
&\quad\times \left( \sum_{j=0}^k \frac{r^j t^j}{j!} \sum_{n_2=k}^\infty S_2(n_2-j,k-j) \frac{t^{n_2-j}}{(n_2-j)!} \right)
\end{split}\end{equation*}
\begin{equation}\begin{split}\label{26}
&= \left( \sum_{n_1=m}^\infty \left( \sum_{l=0}^m {n_1 \choose l} r^l S_2(n_1-l,m-l) \right) \frac{t^{n_1}}{n_1!} \right)\\
&\quad \times \left( \sum_{n_2=k}^\infty \left( \sum_{j=0}^k {n_2 \choose j} r^j S_2(n_2-j,k-j) \right) \frac{t^{n_2}}{n_2!} \right)\\
&=\sum_{n=m+k}^\infty \Bigg\{ \sum_{n_1=m}^n \sum_{l=0}^m \sum_{j=0}^k {n_1 \choose l} {n-n_1 \choose j} r^{l+j} {n \choose n_1} \\
&\qquad \times S_2(n_1-l,m-l) S_2(n-n_1-j,k-j) \Bigg\}\frac{t^n}{n!}.
\end{split}\end{equation}
By \eqref{24} and \eqref{26}, we get
\begin{equation}\begin{split}\label{27}
&{m+k \choose m} S_{2,r} (n,m+k)\\
&=\sum_{n_1=m}^n \sum_{l=0}^m \sum_{j=0}^k {n_1 \choose l} {n-n_1 \choose j}{n \choose n_1} r^{l+j}  S_2(n_1-l,m-l) S_2(n-n_1-j,k-j).
\end{split}\end{equation}
With $r=0$ in \eqref{24}, we have
\begin{equation}\begin{split}\label{28}
\frac{1}{m!}(e^t-1)^m\frac{1}{k!}(e^t-1)^k = {k+m \choose m} \sum_{n=k+m}^\infty S_2(n,m+k) \frac{t^n}{n!},
\end{split}\end{equation}
where $n,m,k \geq 0$. From \eqref{02} and \eqref{09}, we have
\begin{equation*}\begin{split}
&\frac{1}{m!} (e^t-1)^m \frac{1}{k!}(e^t-1)^k = \frac{1}{m!} (e^t-1+rt-rt)^m \frac{1}{k!} (e^t-1+rt-rt)^k \\
&=\left( \frac{1}{m!} \sum_{l=0}^m {m \choose l } (e^t-1+rt)^{m-l} (-rt)^l \right)
\left( \frac{1}{k!} \sum_{j=0}^k {k \choose j } (e^t-1+rt)^{k-j} (-rt)^j \right) \\
&= \left( \sum_{l=0}^m \frac{(-1)^l r^l}{l!} t^l  \sum_{n_1=m}^\infty S_{2,r}(n_1-l,m-l) \frac{t^{n_1-l}}{(n_1-l)!}\right) \\
&\quad \times \left( \sum_{j=0}^k \frac{(-1)^jr^j}{j!}t^j \sum_{n_2=k}^\infty  S_{2,r}(n_2-j,k-j) \frac{t^{n_2-j}}{(n_2-j)!}\right)
\end{split}\end{equation*}
\begin{equation}\begin{split}\label{29}
&= \left( \sum_{n_1=m}^\infty \sum_{l=0}^m {n_1 \choose l} (-1)^l r^l S_{2,r}(n_1-l,m-l) \frac{t^{n_1}}{n_1!}   \right)\\
&\quad \times  \left( \sum_{n_2=k}^\infty \sum_{j=0}^k {n_2 \choose j} (-1)^j r^j S_{2,r}(n_2-j,k-j) \frac{t^{n_2}}{n_2!}   \right)\\
&= \sum_{n=m+k}^\infty \Bigg\{ \sum_{n_1=m}^n \sum_{l=0}^m \sum_{j=0}^k {n_1 \choose l} {n-n_1 \choose j} {n \choose n_1} (-1)^{l+j} r^{l+j}\\
&\qquad	\times  S_{2,r}(n_1-l,m-l) S_{2,r}(n-n_1-j,k-j) \Bigg\} \frac{t^n}{n!}.
\end{split}\end{equation}
Comparing the coefficients on both sides of \eqref{28} and \eqref{29}, we have
\begin{equation}\begin{split}\label{30}
&{k+m \choose m}  S_2(n,m+k) \\
&=\sum_{n_1=m}^n \sum_{l=0}^m \sum_{j=0}^k {n_1 \choose l} {n-n_1 \choose j} {n \choose n_1} (-1)^{l+j} r^{l+j} \\
&\quad \times S_{2,r}(n_1-l,m-l) S_{2,r}(n-n_1-j,k-j),
\end{split}\end{equation}
where $n,m,k \geq 0$ with $n \geq m+k$.

\section{Further Remarks}

A random variable $X$, taking on one of the values $0,1,2,\cdots,$ is said to be a Poisson random variable with parameter $\lambda >0$ if $P(i)=P(X=i)=e^{-\lambda} \frac{\lambda^i}{i!}$, $i=0,1,2,\cdots$. Note that $\sum_{i=0}^\infty P(i) = e^{-\lambda} \sum_{i=0}^\infty \frac{\lambda^i}{i!} = e^{-\lambda}e^\lambda =1$.

The Bell polynomials $Bel_n(x)$, $(n \geq 0)$, are known to be connected with the Poisson distribution. More precisely, $Bel_n(\lambda)$ can be expressed in terms of the moments of Poisson random variable $x$ with parameter $\lambda >0$ as
\begin{equation*}\begin{split}
Bel_n(\lambda) = E[X^n],\,\,(n \in \mathbb{N}).
\end{split}\end{equation*}

Let $X$ be a Poisson random variable with paramerer $\lambda>0$. Then we observe that

\begin{equation}\begin{split}\label{31}
&E[e^{t(X+r\lambda)}] = \sum_{n=0}^\infty E[(X+r\lambda)^n ] \frac{t^n}{n!}\\
&= \sum_{n=0}^\infty \left( \sum_{x=0}^\infty (x+r\lambda)^n \frac{\lambda^x}{x!} e^{-\lambda} \right) \frac{t^n}{n!} \\
&= e^{-\lambda} \sum_{x=0}^\infty \left( \sum_{n=0}^\infty (x+r\lambda)^n \frac{t^n}{n!} \right) \frac{\lambda^x}{x!} \\
&= e^{-\lambda} \sum_{x=0}^\infty e^{(x+r\lambda)t} \frac{\lambda^x}{x!} = e^{rt\lambda-\lambda} \sum_{x=0}^\infty e^{xt} \frac{\lambda^x}{x!}\\
&= e^{\lambda(e^t-1+rt)} = \sum_{n=0}^\infty Bel_{n,r}(\lambda) \frac{t^n}{n!}.
\end{split}\end{equation}
Thus, by \eqref{31}, we see that the extended Bell polynomials are expressed in terms of the moments of Poisson random variable $X$ with parameter $\lambda>0$ as follows:

\begin{equation}\begin{split}\label{32}
E[(X+r\lambda)^n ] = Bel_{n,r}(\lambda),
\end{split}\end{equation}
where $n \in \mathbb{N}$ and $r \in \mathbb{R}$. By binomial theorem, we get
\begin{equation}\begin{split}\label{33}
(X+r\lambda)^n = \sum_{l=0}^n {n \choose l} r^l \lambda^l X^{n-l}.
\end{split}\end{equation}
Thus, by \eqref{33}, we get
\begin{equation}\begin{split}\label{34}
Bel_{n,r}(\lambda) &= E[(X+r\lambda)^n]=\sum_{l=0}^n {n \choose l} r^l \lambda^l E[X^{n-l}]\\
&= \sum_{l=0}^n {n \choose l} r^l \lambda^l Bel_{n-l}(\lambda).
\end{split}\end{equation}
From \eqref{09} and \eqref{31}, we note that
\begin{equation}\begin{split}\label{35}
&\sum_{n=0}^\infty E[(X+r\lambda)^n] \frac{t^n}{n!} = e^{\lambda(e^t-1+rt)}\\
&= \sum_{m=0}^\infty \lambda^m \frac{1}{m!} (e^t-1+rt)^m = \sum_{m=0}^\infty \lambda^m \sum_{n=m}^\infty S_{2,r}(n,m) \frac{t^n}{n!} \\
&= \sum_{n=0}^\infty \left( \sum_{m=0}^n \lambda^m S_{2,r}(n,m) \right) \frac{t^n}{n!}.
\end{split}\end{equation}

Thus, by comparing the coefficients on both sides of \eqref{35}, we get
\begin{equation}\begin{split}\label{36}
E[(X+r\lambda)^n ] = \sum_{m=0}^n \lambda^m S_{2,r}(n,m)=Bel_{n,r}(\lambda),
\end{split}\end{equation}
where $n \in \mathbb{N} \cup \{0\}$ and $X$ is a Poisson random variable with parameter $\lambda>0$.

Now, we observe that
\begin{equation}\begin{split}\label{37}
&e^{tx}E[e^{t(X+r\lambda)}]\\
&=\left( \sum_{l=0}^\infty \frac{x^l}{l!} t^l \right) \left( \sum_{m=0}^\infty E[(X+r\lambda)^m ] \frac{t^m}{m!} \right)\\
&= \sum_{n=0}^\infty \left( \sum_{m=0}^n {n \choose m} x^{n-m} E[(X+r\lambda)^m ] \right) \frac{t^n}{n!}.
\end{split}\end{equation}
On the other hand,
\begin{equation}\begin{split}\label{38}
&e^{tx}E[e^{t(X+r\lambda)}] = e^{\lambda(e^t-1+rt)}e^{xt}\\
&= \sum_{k=0}^\infty \lambda^k \frac{1}{k!} (e^t-1+rt)^k e^{xt}\\
&= \sum_{k=0}^\infty \lambda^k \sum_{n=k}^\infty S_{2,r}(n,k|x) \frac{t^n}{n!}
= \sum_{n=0}^\infty \left( \sum_{k=0}^n \lambda^k S_{2,r}(n,k|x) \right) \frac{t^n}{n!}.
\end{split}\end{equation}
Thus, by \eqref{37} and \eqref{38}, we get
\begin{equation}\begin{split}\label{39}
&\sum_{m=0}^n {n \choose m} x^{n-m} E[(X+r\lambda)^m]= \sum_{k=0}^n \lambda^k S_{2,r}(n,k|x), 
\end{split}\end{equation}
where $n,k \geq 0$ and $X$ is Poisson random variable with parameter $\lambda>0$.

The \eqref{39} is equivalent to
\begin{equation}\begin{split}\label{40}
&\sum_{m=0}^n {n \choose m} x^{n-m} Bel_{m,r}(\lambda)\\
&=\sum_{k=0}^n \lambda^k S_{2,r} (n,k|x), \,\,\text{where}\,\,n \geq 0, \, r \in \mathbb{R}.
\end{split}\end{equation}

\end{document}